\documentclass[a4paper]{jpconf} 
\usepackage{iopams,amsgen,amsfonts,amssymb,amsbsy}
\input cyracc.def
\font\tencyr=wncyr10
\font\tencyi=wncyi10
\def\cyr{\tencyr\cyracc}
\def\cyi{\tencyi\cyracc}
\def\mysavedown#1{\edef\mysubs{\mysubs#1}}
\def\mysaveup#1{\edef\mysups{\mysups#1}}
\def\mydown#1{{\mytensor}_{\vphantom{\mysubs}#1}}
\def\myup#1{{\mytensor}^{\vphantom{\mysups}#1}}
\def\tensor#1#2{
  #1
  \def\mytensor{\vphantom{#1}}
  \def\mysubs{\relax}
  \def\mysups{\relax}
  \let\down=\mysavedown
  \let\up=\mysaveup
  #2
  \let\down=\mydown
  \let\up=\myup
  #2
  }

\newtheorem{Lem}{Lemma}[section] 
\newtheorem{Prop}[Lem]{Proposition} 
\newtheorem{Thm}[Lem]{Theorem} 

\newtheorem{Def}[Lem]{Definition}

\newenvironment{proof}{{\sc Proof.}}{$\Box$}

\newcommand{\bbN}{{\mathbb N}} 
\newcommand{\bbK}{{\mathbb K}} 
\newcommand{\calA}{{\mathcal A}} 
 
\newcommand{\calH}{{\mathcal H}} 
\newcommand{\calL}{{\mathcal L}}

\newcommand{\calS}{{\mathcal S}} 
\newcommand{\calT}{{\mathcal T}} 
\newcommand{\calV}{{\mathcal V}} 
 
\newcommand{\id}{\mathrm{id}} 
\newcommand{\frakl}{{\mathfrak l}} 
\newcommand{\frakr}{{\mathfrak r}} 

\newcommand{\fR}{{\mathfrak R}}

\newcommand{\regrep}{\breve\omega}

\begin{document} 

\title{Stationary or static space-times and Young tableaux} 
\author{Bernd Fiedler}
\address{Eichelbaumstr. 13, D-04249 Leipzig, Germany. URL: http://www.fiemath.de/}  
\ead{bfiedler@fiemath.de}  


\section{Introduction}
The subject of the present paper is the search for examples of Riemann tensors which contain terms with a structure of certain generators of algebraic curvature tensors.

Let $V$ be a finite-dimensional $\bbK$-vector space, where $\bbK$ is the field of real or complex numbers. We denote by $\calT_r V$ the $\bbK$-vector space of covariant tensors of order $r$ over $V$.
\begin{Def}
The $\bbK$-vector space $\calA(V)\subset\calT_4 V$ of all \textit{algebraic curvature tensors} is the set of all tensors $\fR\in\calT_4 V$ which satisfy for all $w, x, y, z\in V$
\begin{eqnarray}
 & & \fR(w,x,y,z)\;=\; - \fR(w,x,z,y)\;=\;\fR(y,z,w,x) \\
 & & \fR(w,x,y,z) + \fR(w,y,z,x) + \fR(w,z,x,y)\;=\; 0\,.
\end{eqnarray}
\end{Def}
For algebraic curvature tensors several types of generators are known. For instance, algebraic curvature tensors can be generated by the following tensors:
\begin{eqnarray}
\gamma(S)_{\kappa\lambda\mu\nu} & := & \label{e1.3}%
S_{\kappa\nu}S_{\lambda\mu} - S_{\kappa\mu}S_{\lambda\nu}\hspace{3cm},\;S\in\calS^2(V)\,, \\
\alpha(A)_{\kappa\lambda\mu\nu} & := & \label{e1.4}%
2 A_{\kappa\lambda}A_{\mu\nu} + A_{\kappa\mu}A_{\lambda\nu} - A_{\kappa\nu}A_{\lambda\mu}\hspace{0.72cm},\;A\in{\Lambda}^2(V)\,,
\end{eqnarray}
where $\calS^p(V)$, ${\Lambda}^p(V)$ denotes the spaces of totally symmetric/alternating $p$-forms over $V$. P. Gilkey \cite[pp.41-44, P.236]{gilkey5} and B. Fiedler \cite{fie20} gave different proofs for
\begin{Thm}
$\calA(V) = \mathrm{Span}_{S\in\calS^2(V)}\{\gamma(S)\} = \mathrm{Span}_{A\in{\Lambda}^2(V)}\{\alpha(A)\}$.
\end{Thm}
But $\calA(V)$ possesses also generators on the basis of products $U\otimes w$ or $w\otimes U$, $U\in\calT_3 V$, $w\in\calT_1 V$, where $U$ has a so-called irreducible $(2\,1)$-symmetry. B. Fiedler proved in \cite{fie04a}
\begin{Thm}\footnote{See section \ref{sec4} for details about symmetry classes of tensors, Young tableaux and Young symmetrizers.} \label{thm1.3}%
Let $\frakr\subset\bbK[\calS_3]$ be a minimal right ideal belonging to the partition $(2\,1)\vdash 3$ and let $\calT_{\frakr}$ be the symmetry class of tensors $U\in\calT_3 V$ that is defined by $\frakr$. Further, let $y_t$ be the Young symmetrizer of the Young tableau
\begin{eqnarray} \label{e4.8}%
t & := &
\begin{array}{|c|c|}
\hline
1 & 3 \\
\hline
2 & 4 \\
\hline
\end{array}\,.
\end{eqnarray}
Then the following statements are equivalent
\begin{enumerate}
\item{$\calA(V) = \mathrm{Span}_{U\in\calT_{\frakr} , w\in\calT_1 V}\{y_t^{\ast}(U\otimes w)\} = \mathrm{Span}_{U\in\calT_{\frakr} , w\in\calT_1 V}\{y_t^{\ast}(w\otimes U)\}$.}
\item{$\frakr$ is different from the right ideal
$\frakr_0 := f_0\cdot\bbK[\calS_3]$ which is generated by the idempotent
$f_0 := \frac{1}{2}\left\{\id - (1\,3)\right\} -
\frac{1}{6}\,\sum_{p\in\calS_3}\,\mathrm{sign}(p)\,p$.
}
\end{enumerate}
\end{Thm}
In differential geometry or general relativity theory (GR), many examples of Riemann tensors are known, in which expressions of the type (\ref{e1.3}) or (\ref{e1.4}) occur (see e.g. \cite[Thm.2.1]{gilkey3},\cite[Thm.2]{drgr04},\cite[Sec.3]{dfgg04}). In the present paper we want to show that there exist also curvature formulas in which generators from Theorem \ref{thm1.3} appear realized by differentiable tensor fields in a natural way. We found such generators in curvature formulas of stationary and static space-times of GR. 

Let $(M,g)$ be a 4-dimensional pseudo-Riemannian manifold of class $C^{\infty}$ whose fundamental tensor $g$ has signature $(+---)$. Let $\nabla$ be the Levi-Civita connection of $g$.
\begin{Def}
\begin{enumerate}
\item{$(M,g)$ is called a \textit{stationary} space-time if there exists a timelike Killing field $\xi$ on $M$ characterized by the conditions
\begin{eqnarray} \label{e1.1}%
{\xi}_{\mu ; \nu} + {\xi}_{\nu ; \mu}\;=\; 0 & \;\;,\;\; &
g_{\mu\nu} {\xi}^{\mu} {\xi}^{\nu}\;>\; 0\,.
\end{eqnarray}}
\item{A stationary $(M,g)$ is called a \textit{static} space-time if the Killing field $\xi$ is \textit{hypersurface-orthogonal}, i.e. $\xi$ satisfies the additional condition
${\xi}_{[\mu;\nu} {\xi}_{\lambda]} = 0$.
}
\end{enumerate}
\end{Def}
If $(M,g)$ is stationary then one can construct local coordinates $t, x^1, x^2, x^3$ around every point $p\in M$ such that
$g_{\mu\nu} = g_{\mu\nu}(x^1, x^2, x^3)$ and $\xi\;=\;{\partial}_t$.
If $(M,g)$ is static, then we can choose these local coordinates $t, x^1, x^2, x^3$ in such a way that
\begin{eqnarray}
ds^2\;=\;g_{\mu\nu}dx^{\mu}dx^{\nu} & = & f(x^1, x^2, x^3)^2dt^2 - h_{ab}(x^1, x^2, x^3)dx^adx^b \\
\mu , \nu = 0 , 1 , 2 , 3
& &
a , b = 1 , 2 , 3\,. \nonumber
\end{eqnarray}
Here $d{\sigma}^2 := h_{ab}dx^adx^b$ is a positive definite, 3-dimensional metric.

Now we formulate the main results of our paper.
\begin{Prop} \label{prop1.2}%
Let ${\tau}^{\mu} := ({\xi}_{\alpha}{\xi}^{\alpha})^{-1/2}{\xi}^{\mu}$ be the timelike unit vector field which is proportional to the above Killing field ${\xi}^{\mu}$. Further, let $Z_{\kappa\lambda\mu\nu}$ be the {\rm (}covariant{\rm )} spatial projection {\rm (\ref{e2.26})} of the covariant Riemann tensor $R_{\kappa\lambda\mu\nu}$, let $P_{\kappa\lambda\mu\nu}$ be the {\rm (}covariant{\rm ) 3}-dimensional curvature tensor from {\rm (\ref{e2.24})} and let $F_{\mu}$ be the field {\rm (\ref{e2.13})}. Then the following formula holds:
\begin{eqnarray}
- Z_{\kappa\lambda\mu\nu}\;=\;P_{\kappa\lambda\mu\nu} & - & 
  \frac{3}{4} \tensor{F}{\down{\lambda}} 
    \tensor{F}{\down{\nu}} \tensor{\tau}{\down{\kappa}} \tensor{\tau}{\down{\mu}} + 
  \frac{3}{4} \tensor{F}{\down{\kappa}} 
    \tensor{F}{\down{\nu}} \tensor{\tau}{\down{\lambda}} \tensor{\tau}{\down{\mu}} + 
  \frac{3}{4} \tensor{F}{\down{\lambda}} 
    \tensor{F}{\down{\mu}} \tensor{\tau}{\down{\kappa}} \tensor{\tau}{\down{\nu}} - 
  \frac{3}{4} \tensor{F}{\down{\kappa}} 
    \tensor{F}{\down{\mu}} \tensor{\tau}{\down{\lambda}} \tensor{\tau}{\down{\nu}} - \nonumber \\ & &
  \tensor{F}{\down{\nu}} \tensor{\tau}{\down{\mu}} 
    \tensor{\tau}{\down{[\kappa;}\down{\lambda]}} + 
  \tensor{F}{\down{\mu}} \tensor{\tau}{\down{\nu}} 
    \tensor{\tau}{\down{[\kappa;}\down{\lambda]}} - 
  \frac{1}{2} \tensor{F}{\down{\nu}} 
    \tensor{\tau}{\down{\lambda}} \tensor{\tau}{\down{[\kappa;}\down{\mu]}} + 
  \frac{1}{2} \tensor{F}{\down{\lambda}} 
    \tensor{\tau}{\down{\nu}} \tensor{\tau}{\down{[\kappa;}\down{\mu]}} + \nonumber \\ & &
  \frac{1}{2} \tensor{F}{\down{\mu}} 
    \tensor{\tau}{\down{\lambda}} \tensor{\tau}{\down{[\kappa;}\down{\nu]}} - 
  \frac{1}{2} \tensor{F}{\down{\lambda}} 
    \tensor{\tau}{\down{\mu}} \tensor{\tau}{\down{[\kappa;}\down{\nu]}} + 
  \frac{1}{2} \tensor{F}{\down{\nu}} 
    \tensor{\tau}{\down{\kappa}} \tensor{\tau}{\down{[\lambda;}\down{\mu]}} - 
  \frac{1}{2} \tensor{F}{\down{\kappa}} 
    \tensor{\tau}{\down{\nu}} \tensor{\tau}{\down{[\lambda;}\down{\mu]}} - \label{e1.5} \\ & &
  \frac{1}{2} \tensor{F}{\down{\mu}} 
    \tensor{\tau}{\down{\kappa}} \tensor{\tau}{\down{[\lambda;}\down{\nu]}} + 
  \frac{1}{2} \tensor{F}{\down{\kappa}} 
    \tensor{\tau}{\down{\mu}} \tensor{\tau}{\down{[\lambda;}\down{\nu]}} - 
  \tensor{F}{\down{\lambda}} \tensor{\tau}{\down{\kappa}} 
    \tensor{\tau}{\down{[\mu;}\down{\nu]}} + 
  \tensor{F}{\down{\kappa}} \tensor{\tau}{\down{\lambda}} 
    \tensor{\tau}{\down{[\mu;}\down{\nu]}} - \nonumber \\ & &
  2 \tensor{\tau}{\down{[\kappa;}\down{\lambda]}} \tensor{\tau}{\down{[\mu;}\down{\nu]}} + 
  \tensor{\tau}{\down{[\kappa;}\down{\nu]}} \tensor{\tau}{\down{[\lambda;}\down{\mu]}} - 
  \tensor{\tau}{\down{[\kappa;}\down{\mu]}} \tensor{\tau}{\down{[\lambda;}\down{\nu]}} \nonumber
\end{eqnarray}
\end{Prop}
\vspace{0.3cm}
\begin{Thm} \label{thm1.6}%
Let $M_p$ be the tangent space of $M$ in a point $p\in M$.
\begin{enumerate}
\item{In every $p\in M$ the tensor $({\tau}_{\lambda}{\tau}_{[\mu ; \nu]})|_p$ belongs to a symmetry class $\calT_{\frakr}\subseteq\calT_3 M_p$ whose defining right ideal $\frakr\subset\bbK[\calS_3]$ has a decomposition $\frakr = \frakr_1\oplus\frakr_2$ into {\rm 2} minimal right ideals $\frakr_i$ which is described by the Littlewood-Richardson product
\begin{eqnarray}
[1^2][1] & \sim & [2\,1] + [1^3]\,. \label{e1.6}
\end{eqnarray}}
\item{In every $p\in M$ the tensor $(F_{\kappa}{\tau}_{\lambda}{\tau}_{[\mu ; \nu]})|_p$ lies in a symmetry class $\calT_{\frakr}\subseteq\calT_4 M_p$ whose defining right ideal $\frakr\subset\bbK[\calS_4]$ has a decomposition $\frakr = \frakr_1\oplus\ldots\oplus\frakr_5$ into {\rm 5} minimal right ideals $\frakr_i$ which is described by the Littlewood-Richardson products
\begin{eqnarray}
[2\,1][1]\;\sim\;[3\,1] + [2^2] + [2\,1^2]
& \;\;,\;\; &
[1^3][1]\;\sim\;[2\,1^2] + [1^4]\,. \label{e1.14}%
\end{eqnarray}
At most the product $[2\,1][1]$ yields a contribution to the symmetry class $\calA(V)$ of algebraic curvature tensors which belongs to $[2^2]$.}
\end{enumerate}
\end{Thm}
\begin{Lem} \label{lemma1.4}%
For the above fields ${\xi}^{\mu}$ and ${\tau}^{\mu}$ the conditions
${\xi}_{[\lambda}{\xi}_{\mu ; \nu]} = 0$ and
${\tau}_{[\lambda}{\tau}_{\mu ; \nu]} = 0$
are equivalent.
\end{Lem}
\begin{Thm} \label{thm1.8}%
A stationary space-time is static iff ${\tau}_{[\lambda}{\tau}_{\mu ; \nu]} = 0$ or, equivalently, iff the part $[1^3]$ of {\rm (\ref{e1.6})} vanishes. In a static space-time the tensor ${\tau}_{\lambda}{\tau}_{[\mu ; \nu]}$ belongs to an irreducible $(2\,1)$-symmetry class.
\end{Thm}
Because of (\ref{e1.6}) the tensor ${\tau}_{\lambda}{\tau}_{[\mu ; \nu]}$ possesses a decomposition
\begin{eqnarray} \label{e1.17}%
{\tau}_{\lambda}{\tau}_{[\mu ; \nu]} & = &
{\theta}_{\lambda\mu\nu} + {\tau}_{[\lambda}{\tau}_{\mu ; \nu]}\,,
\end{eqnarray}
where ${\theta}_{\lambda\mu\nu} := {\tau}_{\lambda}{\tau}_{[\mu ; \nu]} -
{\tau}_{[\lambda}{\tau}_{\mu ; \nu]}$ is the unique part of ${\tau}_{\lambda}{\tau}_{[\mu ; \nu]}$ which has an irreducible $(2\,1)$-symmetry. (\ref{e1.14}) leads to the remarkable consequence that ${\tau}_{[\lambda}{\tau}_{\mu ; \nu]}$ does not yield a contribution to (\ref{e1.5}) even in the case of a stationary space-time.
\begin{Thm} \label{thm1.9}%
The substitution {\rm (\ref{e1.17})} transforms {\rm (\ref{e1.5})} into
\begin{eqnarray}
- Z_{\kappa\lambda\mu\nu}\;=\;P_{\kappa\lambda\mu\nu} & - &
  \frac{3}{4} \tensor{F}{\down{\lambda}} 
    \tensor{F}{\down{\nu}} \tensor{\tau}{\down{\kappa}} \tensor{\tau}{\down{\mu}} + 
  \frac{3}{4} \tensor{F}{\down{\kappa}} 
    \tensor{F}{\down{\nu}} \tensor{\tau}{\down{\lambda}} \tensor{\tau}{\down{\mu}} + 
  \frac{3}{4} \tensor{F}{\down{\lambda}} 
    \tensor{F}{\down{\mu}} \tensor{\tau}{\down{\kappa}} \tensor{\tau}{\down{\nu}} - 
  \frac{3}{4} \tensor{F}{\down{\kappa}} 
    \tensor{F}{\down{\mu}} \tensor{\tau}{\down{\lambda}} \tensor{\tau}{\down{\nu}} - \nonumber \\ & &
  \tensor{F}{\down{\nu}} \tensor{\theta}{\down{\mu}\down{\kappa}\down{\lambda}} + 
  \tensor{F}{\down{\mu}} \tensor{\theta}{\down{\nu}\down{\kappa}\down{\lambda}} - 
  \frac{1}{2} \tensor{F}{\down{\nu}} 
    \tensor{\theta}{\down{\lambda}\down{\kappa}\down{\mu}} + 
  \frac{1}{2} \tensor{F}{\down{\lambda}} 
    \tensor{\theta}{\down{\nu}\down{\kappa}\down{\mu}} + \nonumber \\ & &
  \frac{1}{2} \tensor{F}{\down{\mu}} 
    \tensor{\theta}{\down{\lambda}\down{\kappa}\down{\nu}} - 
  \frac{1}{2} \tensor{F}{\down{\lambda}} 
    \tensor{\theta}{\down{\mu}\down{\kappa}\down{\nu}} + 
  \frac{1}{2} \tensor{F}{\down{\nu}} 
    \tensor{\theta}{\down{\kappa}\down{\lambda}\down{\mu}} - 
  \frac{1}{2} \tensor{F}{\down{\kappa}} 
    \tensor{\theta}{\down{\nu}\down{\lambda}\down{\mu}} - \label{e1.18} \\ & &
  \frac{1}{2} \tensor{F}{\down{\mu}} 
    \tensor{\theta}{\down{\kappa}\down{\lambda}\down{\nu}} + 
  \frac{1}{2} \tensor{F}{\down{\kappa}} 
    \tensor{\theta}{\down{\mu}\down{\lambda}\down{\nu}} - 
  \tensor{F}{\down{\lambda}} \tensor{\theta}{\down{\kappa}\down{\mu}\down{\nu}} + 
  \tensor{F}{\down{\kappa}} \tensor{\theta}{\down{\lambda}\down{\mu}\down{\nu}} - \nonumber \\ & &
  2 \tensor{\tau}{\down{[\kappa;}\down{\lambda]}} \tensor{\tau}{\down{[\mu;}\down{\nu]}} + 
  \tensor{\tau}{\down{[\kappa;}\down{\nu]}} \tensor{\tau}{\down{[\lambda;}\down{\mu]}} - 
  \tensor{\tau}{\down{[\kappa;}\down{\mu]}} \tensor{\tau}{\down{[\lambda;}\down{\nu]}} \nonumber
\end{eqnarray}
In {\rm (\ref{e1.18})} the tensor ${\tau}_{[\lambda}{\tau}_{\mu ; \nu]}$ does not appear.
\end{Thm}
\begin{Thm} \label{thm1.10}%
Let ${\theta}_{\lambda\mu\nu} := {\tau}_{\lambda}{\tau}_{[\mu ; \nu]} - {\tau}_{[\lambda}{\tau}_{\mu ; \nu]}$ be the tensor field from {\rm (\ref{e1.17})}.
\begin{enumerate}
\item{The symmetry properties of $\theta$ are described by the relations
\begin{eqnarray}
 & & \label{e1.19}%
\begin{array}{cccccccl}
0 & = & & \tensor{\theta}{\down{\lambda}\down{\mu}\down{\nu}} & & & + & \hspace*{2.1cm}\tensor{\theta}{\down{\lambda}\down{\nu}\down{\mu}} \\
0 & = & - & \tensor{\theta}{\down{\lambda}\down{\mu}\down{\nu}} & + &
  \tensor{\theta}{\down{\nu}\down{\mu}\down{\lambda}} & + & \hspace*{1.4cm}\tensor{\theta}{\down{\mu}\down{\lambda}\down{\nu}} \\
0 & = & & \tensor{\theta}{\down{\lambda}\down{\mu}\down{\nu}} & - & \tensor{\theta}{\down{\nu}\down{\mu}\down{\lambda}} & + & \hspace*{0.7cm}\tensor{\theta}{\down{\mu}\down{\nu}\down{\lambda}} \\
0 & = & & & & \tensor{\theta}{\down{\nu}\down{\mu}\down{\lambda}} & + & \tensor{\theta}{\down{\nu}\down{\lambda}\down{\mu}} \\
\end{array}
\end{eqnarray}
}
\item{Applying {\rm (\ref{e1.19})} we can reduce the {\rm 20} summands of {\rm (\ref{e1.18})} to the {\rm 12} summands
\begin{eqnarray}
- Z_{\kappa\lambda\mu\nu}\;=\;P_{\kappa\lambda\mu\nu} & - &
\frac{3}{4}  \tensor{F}{\down{\lambda}} 
    \tensor{F}{\down{\nu}} \tensor{\tau}{\down{\kappa}} \tensor{\tau}{\down{\mu}} + 
  \frac{3}{4} \tensor{F}{\down{\kappa}} 
    \tensor{F}{\down{\nu}} \tensor{\tau}{\down{\lambda}} \tensor{\tau}{\down{\mu}} + 
  \frac{3}{4} \tensor{F}{\down{\lambda}} 
    \tensor{F}{\down{\mu}} \tensor{\tau}{\down{\kappa}} \tensor{\tau}{\down{\nu}} - 
  \frac{3}{4} \tensor{F}{\down{\kappa}} 
    \tensor{F}{\down{\mu}} \tensor{\tau}{\down{\lambda}} \tensor{\tau}{\down{\nu}} \nonumber \\ & &
+  \tensor{\tau}{\down{[\kappa;}\down{\nu]}} \tensor{\tau}{\down{[\lambda;}\down{\mu]}} - 
  \tensor{\tau}{\down{[\kappa;}\down{\mu]}} \tensor{\tau}{\down{[\lambda;}\down{\nu]}} - 
  2 \tensor{\tau}{\down{[\kappa;}\down{\lambda]}} \tensor{\tau}{\down{[\mu;}\down{\nu]}} \label{e1.20} \\ & &
\hspace*{-1.4cm}\underbrace{
-  \frac{3}{2} \tensor{F}{\down{\lambda}} \tensor{\theta}{\down{\kappa}\down{\mu}\down{\nu}} + 
  \frac{3}{2} \tensor{F}{\down{\kappa}} \tensor{\theta}{\down{\lambda}\down{\mu}\down{\nu}} + 
  \frac{3}{2} \tensor{F}{\down{\nu}} \tensor{\theta}{\down{\mu}\down{\lambda}\down{\kappa}} - 
  \frac{3}{2} \tensor{F}{\down{\mu}} \tensor{\theta}{\down{\nu}\down{\lambda}\down{\kappa}}
}_{\hspace*{5cm}= \frac{1}{2} y_t^{\ast}(\theta\otimes F)_{\kappa\lambda\mu\nu}\;\mathrm{according}\;\mathrm{to}\;\mathrm{Theorem}\;\ref{thm1.3}}\hspace{-1cm}. \nonumber
\end{eqnarray}
}
\end{enumerate}
\end{Thm}

\section{A projection formalism} \label{sec2}
In investigations of stationary or static space-times one can use a projection formalism which is described for instance in \cite[pp. 180]{kramstepcalherlt} or \cite[pp.49]{vlad}. We apply the formalism from the book \cite[pp.49]{vlad} by Yu. S. Vladimirov in our paper.

The formalism of Vladimirov starts with the assumtion that a timelike unit vector field ${\tau}^{\mu}$ is given on $M$ which describes the 4-speed of a ''\textit{continuum of observers}''. If we define
\begin{eqnarray} \label{e2.1}%
h_{\mu\nu}\;:=\;{\tau}_{\mu}{\tau}_{\nu} - g_{\mu\nu} & , &
h^{\mu}_{\nu}\;:=\; g^{\mu\alpha}h_{\alpha\nu}\;\;\;,\;\;\;
h^{\mu\nu}\;:=\; g^{\mu\alpha}g^{\nu\beta}h_{\alpha\beta}\,,
\end{eqnarray}
then we obtain the following decompositions of the metric tensors:
\begin{eqnarray} \label{e2.2}%
g_{\mu\nu}\;:=\;{\tau}_{\mu}{\tau}_{\nu} - h_{\mu\nu} & , &
g^{\mu\nu}\;:=\;{\tau}^{\mu}{\tau}^{\nu} - h^{\mu\nu}\,.
\end{eqnarray}
A simple consequence of
${\tau}^{\mu}{\tau}_{\mu} = g_{\mu\nu}{\tau}^{\mu}{\tau}^{\nu} = 1$ is
\begin{eqnarray} \label{e2.3}%
h_{\mu\nu}{\tau}^{\nu}\;=\;0 & , &
h^{\mu}_{\nu}{\tau}^{\nu}\;=\;0\;\;\;,\;\;\;
h^{\mu\nu}{\tau}_{\nu}\;=\;0\,.
\end{eqnarray}
If $B_{{\mu}_1\ldots{\mu}_r}^{\hspace*{1cm}{\nu}_1\ldots{\nu}_s}$ are the coordinates of a $r$-times \textit{covariant} and $s$-times \textit{contravariant} tensor $B$ then we can define the \textit{time component} of $B$
\begin{eqnarray}
\tilde{B} & := & B_{{\mu}_1\ldots{\mu}_r}^{\hspace*{1cm}{\nu}_1\ldots{\nu}_s}
{\tau}^{{\mu}_1}\ldots{\tau}^{{\mu}_r}{\tau}_{{\nu}_1}\ldots{\tau}_{{\mu}_s}
\end{eqnarray}
and the \textit{spatial projection} of $B$
\begin{eqnarray} \label{e2.5}%
\tilde{B}_{{\alpha}_1\ldots{\alpha}_r}^{\hspace*{1cm}{\beta}_1\ldots{\beta}_s}
& := &
(-1)^{r+s}\,h_{{\alpha}_1}^{{\mu}_1}\ldots h_{{\alpha}_r}^{{\mu}_r}
h_{{\nu}_1}^{{\beta}_1}\ldots h_{{\nu}_s}^{{\beta}_s}\,
B_{{\mu}_1\ldots{\mu}_r}^{\hspace*{1cm}{\nu}_1\ldots{\nu}_s}\,.
\end{eqnarray}
Furthermore we can form mixed projections, for instance
\begin{eqnarray} \label{e2.6}%
\tilde{B}_{{\alpha}_1\ldots{\alpha}_{r-1}}^{\hspace*{1.3cm}{\beta}_1\ldots{\beta}_s}
& := &
(-1)^{r+s-1}\,h_{{\alpha}_1}^{{\mu}_1}\ldots h_{{\alpha}_{r-1}}^{{\mu}_{r-1}}
{\tau}^{{\mu}_r}h_{{\nu}_1}^{{\beta}_1}\ldots h_{{\nu}_s}^{{\beta}_s}\,
B_{{\mu}_1\ldots{\mu}_r}^{\hspace*{1cm}{\nu}_1\ldots{\nu}_s}\,.
\end{eqnarray}
Because of (\ref{e2.3}) the projections (\ref{e2.5}) and (\ref{e2.6}) fulfill
\begin{eqnarray}
{\tau}^{{\alpha}_i}\,\tilde{B}_{\ldots {\alpha}_i\ldots}^{\hspace*{0.9cm}\ldots\ldots}\;=\;0 & \;\;\;,\;\;\; &
{\tau}_{{\beta}_i}\,\tilde{B}_{\ldots\ldots}^{\hspace*{0.6cm}\ldots{{\beta}_i}\ldots}\;=\;0\,.
\end{eqnarray}
If we use the time component
$d\tilde{\tau} := {\tau}_{\mu}dx^{\mu}$
and the spatial projection
$d\tilde{x}^{\nu} := - h_{\mu}^{\nu}dx^{\mu}$ of $dx^{\mu}$ and set
$dl^2 := h_{\mu\nu}d\tilde{x}^{\mu}d\tilde{x}^{\nu}$
then we obtain
\begin{eqnarray}
ds^2 & = & \left({\tau}_{\mu}{\tau}_{\nu} - h_{\mu\nu}\right)dx^{\mu}dx^{\nu}
\;\;=\;\;
d\tilde{\tau}^2 - h_{\mu\nu}d\tilde{x}^{\mu}d\tilde{x}^{\nu}
\;\;=\;\;
d\tilde{\tau}^2 - dl^2\,.
\end{eqnarray}
Finally, (\ref{e2.1}), (\ref{e2.3}) and $4 = g_{\mu\nu}g^{\mu\nu} = 1 + h_{\mu\nu}h^{\mu\nu}$ lead to
\begin{eqnarray}
h_{\alpha}^{\mu}h_{\mu\beta}\;=\;- h_{\alpha\beta} &\;,\;&
h_{\mu\alpha}h^{\alpha\nu}\;=\; - h_{\mu}^{\nu}\;\;\;\;,\;\;\;\;h_{\mu\nu}h^{\mu\nu}\;=\;3\,.
\end{eqnarray}
Now we consider the decomposition
\begin{eqnarray}
{\tau}_{\mu ; \nu} & = & \frac{1}{2}\left({\tau}_{\mu , \nu} - {\tau}_{\nu , \mu}\right) + \frac{1}{2}\left({\tau}_{\mu ; \nu} + {\tau}_{\nu ; \mu}\right)
\end{eqnarray}
of the covariant derivative ${\tau}_{\mu ; \nu}$. If we covariantly differentiate $1 = {\tau}_{\mu}{\tau}^{\mu}$ we obtain immediately
\begin{eqnarray} \label{e2.12}%
{\tau}_{\mu ; \nu}{\tau}^{\mu}\;=\;0 & , &
\left({\tau}_{\mu , \nu} - {\tau}_{\nu , \mu}\right){\tau}^{\mu}{\tau}^{\nu}\;=\;0
\;\;,\;\;
\left({\tau}_{\mu ; \nu} + {\tau}_{\nu ; \mu}\right){\tau}^{\mu}{\tau}^{\nu}\;=\;0\,.
\end{eqnarray}
\begin{Def}
The following quantities $F_{\alpha}, A_{\alpha\beta}, D_{\alpha\beta}$ are usefull in formulas for the Riemann tensor:
\begin{eqnarray} \label{e2.13}%
\hspace*{1.3cm}F_{\alpha} & := &
-\left({\tau}_{\mu , \nu} - {\tau}_{\nu , \mu}\right){\tau}^{\nu}h_{\alpha}^{\mu}
\;=\;
-\left({\tau}_{\mu ; \nu} + {\tau}_{\nu ; \mu}\right){\tau}^{\mu}h_{\alpha}^{\nu}
\;=\;
\left({\tau}_{\alpha , \nu} - {\tau}_{\nu , \alpha}\right){\tau}^{\nu}\,, \\
A_{\alpha\beta} & := & \label{e2.14}%
\frac{1}{2}\left({\tau}_{\mu , \nu} - {\tau}_{\nu , \mu}\right)h_{\alpha}^{\mu}h_{\beta}^{\nu}\;\;\;,\;\;\;
D_{\alpha\beta}\;:=\;
- \frac{1}{2}\left({\tau}_{\mu ; \nu} + {\tau}_{\nu ; \mu}\right)h_{\alpha}^{\mu}h_{\beta}^{\nu}\,.
\end{eqnarray}
\end{Def}
Obviously, $D$ is a symmetric tensor, $D_{\alpha\beta} = D_{\beta\alpha}$, whereas $A$ is skew-symmetric, $A_{\alpha\beta} = - A_{\beta\alpha}$.
\begin{Lem}
The tensor $A_{\alpha\beta}$ satisfies {\rm (}see {\rm\cite[p.51]{vlad})}
\begin{eqnarray} \label{e2.16}%
A_{\alpha\beta} & = & \frac{1}{2}\left({\tau}_{\alpha , \beta} - {\tau}_{\beta , \alpha}\right) + 
\frac{1}{2}\left({\tau}_{\alpha} F_{\beta} - {\tau}_{\beta} F_{\alpha}\right)\,.
\end{eqnarray}
\end{Lem}
\begin{proof}
Taking into account (\ref{e2.1}) and (\ref{e2.14}) we can write
\begin{eqnarray} \label{e2.17}%
A_{\alpha\beta} & = &
\frac{1}{2}\left({\tau}_{\mu , \nu} - {\tau}_{\nu , \mu}\right)
{\tau}^{\mu}{\tau}^{\nu}{\tau}_{\alpha}{\tau}_{\beta}
- \frac{1}{2}\left({\tau}_{\mu , \nu} - {\tau}_{\nu , \mu}\right)
{\tau}^{\mu}{\tau}_{\alpha}{\delta}_{\beta}^{\nu} \\
 & & 
- \frac{1}{2}\left({\tau}_{\mu , \nu} - {\tau}_{\nu , \mu}\right)
{\tau}^{\nu}{\tau}_{\beta}{\delta}_{\alpha}^{\mu}
+ \frac{1}{2}\left({\tau}_{\mu , \nu} - {\tau}_{\nu , \mu}\right)
{\delta}_{\alpha}^{\mu}{\delta}_{\beta}^{\nu}\,, \nonumber
\end{eqnarray}
where ${\delta}_{\mu}^{\nu}$ denotes the Kronecker symbol. The first summand of (\ref{e2.17}) vanishes because of (\ref{e2.12}). But then (\ref{e2.13}) leads to (\ref{e2.16}).
\end{proof}

We use the following definition of the Christoffel symbols and the Riemann tensor of $\nabla$:
\begin{Def}
Let $\nabla$ be the Levi-Civita connection of $g$. Then we define the \textit{Christoffel symbols} ${\Gamma}_{\mu\nu}^{\lambda}$ and the \textit{Riemann tensor} $R_{\mu\nu\alpha}^{\hspace*{0.6cm}\lambda}$ of $\nabla$ by
\begin{eqnarray} \label{e2.18}%
{\Gamma}_{\mu\nu}^{\lambda} & := & \frac{1}{2}\,g^{\lambda\gamma}
\left(
{\partial}_{\mu}g_{\nu\gamma} +
{\partial}_{\nu}g_{\mu\gamma} -
{\partial}_{\gamma}g_{\mu\nu}
\right)\, \\
R_{\mu\nu\alpha}^{\hspace*{0.6cm}\lambda} & := & \label{e2.19}%
{\partial}_{\mu}{\Gamma}_{\nu\alpha}^{\lambda} -
{\partial}_{\nu}{\Gamma}_{\mu\alpha}^{\lambda} +
{\Gamma}_{\nu\alpha}^{\epsilon}{\Gamma}_{\mu\epsilon}^{\lambda} -
{\Gamma}_{\mu\alpha}^{\epsilon}{\Gamma}_{\nu\epsilon}^{\lambda}\,.
\end{eqnarray}
\end{Def}
The book \cite{vlad} by Vladimirov uses (\ref{e2.18}), too, but defines the Riemannian curvature tensor by
\begin{eqnarray} \label{e2.20}%
\breve{R}_{\;\;\alpha\nu\mu}^{\lambda} & := &
{\partial}_{\nu}{\Gamma}_{\alpha\mu}^{\lambda} -
{\partial}_{\mu}{\Gamma}_{\alpha\nu}^{\lambda} +
{\Gamma}_{\alpha\mu}^{\epsilon}{\Gamma}_{\epsilon\nu}^{\lambda} -
{\Gamma}_{\alpha\nu}^{\epsilon}{\Gamma}_{\epsilon\mu}^{\lambda}\,.
\end{eqnarray}
The transformation between (\ref{e2.19}) and (\ref{e2.20}) reads
\begin{eqnarray} \label{e2.21}%
R_{\mu\nu\alpha}^{\hspace*{0.6cm}\lambda} & = &
- \breve{R}_{\;\;\alpha\nu\mu}^{\lambda}\,.
\end{eqnarray}
\begin{Lem}
The Christoffel symbols ${\Gamma}_{\mu\nu}^{\lambda}$ can be expressed by
${\tau}^{\mu}, h_{\mu\nu}, F_{\alpha}, A_{\alpha\beta}$ and $D_{\alpha\beta}$ in the following way
\begin{eqnarray*}
{\Gamma}_{\alpha\beta}^{\mu} & = &
\left(L_{\alpha\beta}^{\mu} + h_{\alpha , \beta}^{\mu}\right) -
{\tau}^{\mu}\left(A_{\alpha\beta} - D_{\alpha\beta} + F_{\alpha}{\tau}_{\beta}\right) + F^{\mu}{\tau}_{\alpha}{\tau}_{\beta}
+ \left({\tau}_{\alpha}A^{\mu}_{\;\;\beta} + {\tau}_{\beta}A^{\mu}_{\;\;\alpha}\right) - {\tau}_{\alpha}{\tau}^{\mu}_{\;\;,\beta}\,,
\nonumber
\end{eqnarray*}
where
$L_{\alpha\beta}^{\mu} :=
\frac{1}{2}\,h^{\mu\epsilon}
\left(
{\partial}_{\beta}h_{\alpha\epsilon} +
{\partial}_{\alpha}h_{\beta\epsilon} -
{\partial}_{\epsilon}h_{\alpha\beta}
\right)$
denotes the ''Christoffel symbols'' of the ''{\rm 3}-dimensional metric'' $h$.
{\rm (}See {\rm\cite[p.53]{vlad}.)}
\end{Lem}
\begin{Def}
We denote by $P$ the tensor
\begin{eqnarray} \label{e2.24}%
P_{\alpha\beta\rho}^{\hspace*{0.6cm}\lambda}\;:=\;
- h_{\alpha}^{\mu}h_{\beta}^{\nu}h_{\rho}^{\sigma}h_{\epsilon}^{\lambda}
\left(
{\partial}_{\mu}\tilde{L}_{\sigma\nu}^{\epsilon} -
{\partial}_{\nu}\tilde{L}_{\sigma\mu}^{\epsilon} +
\tilde{L}_{\kappa\mu}^{\epsilon}\tilde{L}_{\sigma\nu}^{\kappa} -
\tilde{L}_{\kappa\nu}^{\epsilon}\tilde{L}_{\sigma\mu}^{\kappa}
\right)\;\;,\;\;
\tilde{L}_{\sigma\nu}^{\epsilon}\;:=\;
L_{\sigma\nu}^{\epsilon} + h_{\sigma , \nu}^{\epsilon}\,, & &
\end{eqnarray}
which can be considered the curvature tensor assigned to $h$. {\rm (} See {\rm\cite[p.55]{vlad}).}
\end{Def}
The right-hand side of (\ref{e2.24}) is equal to the right-hand side of formula (3.30) in \cite[p.55]{vlad}. We adapted only the left-hand side of (\ref{e2.24}) to our definition (\ref{e2.19}) of the curvature tensor by means of (\ref{e2.21}). Note that $\tilde{L}_{\sigma\nu}^{\epsilon}$ is not symmetric with respect to $\sigma , \nu$ in general.

The Riemann tensor $R$ possesses three spatial projections.
\begin{Def}
We denote by $Z_{\mu\nu\kappa}^{\hspace*{0.6cm}\lambda}$,
$Y_{\mu\nu\kappa}$, $X_{\nu\kappa}$ the following three spatial projections of the Riemann tensor $R$:
\begin{eqnarray} \label{e2.26}%
Z_{\mu\nu\kappa}^{\hspace*{0.6cm}\lambda}\;:=\;
h_{\mu}^{\alpha}h_{\nu}^{\beta}h_{\kappa}^{\gamma}h_{\delta}^{\lambda}
R_{\alpha\beta\gamma}^{\hspace*{0.6cm}\delta}\;\;,\;\;
Y_{\mu\nu\kappa}\;:=\;
h_{\mu}^{\alpha}h_{\nu}^{\beta}h_{\kappa}^{\gamma}{\tau}_{\delta}
R_{\alpha\beta\gamma}^{\hspace*{0.6cm}\delta}\;\;,\;\;
X_{\nu\kappa}\;:=\;
- h_{\nu}^{\beta}h_{\kappa}^{\gamma}{\tau}^{\alpha}{\tau}_{\delta}
R_{\alpha\beta\gamma}^{\hspace*{0.6cm}\delta}\,. & &
\end{eqnarray}
\end{Def}
In the present paper we consider only $Z_{\mu\nu\kappa}^{\hspace*{0.6cm}\lambda}$.
\begin{Prop}
The spatial projection $Z$ of the Riemann tensor $R$ satisfies
{\small
\begin{eqnarray} \label{e2.29}%
- Z_{\epsilon\gamma\kappa}^{\hspace*{0.5cm}\lambda} & = &
P_{\epsilon\gamma\kappa}^{\hspace*{0.5cm}\lambda} +
2\,A^{\lambda}_{\;\;\kappa}A_{\epsilon\gamma}
+
\left(D_{\epsilon}^{\lambda} + A_{\epsilon}^{\;\;\,\lambda}\right)
\left(D_{\gamma\kappa} + A_{\gamma\kappa}\right) -
\left(D_{\gamma}^{\lambda} + A_{\gamma}^{\;\;\,\lambda}\right)
\left(D_{\epsilon\kappa} + A_{\epsilon\kappa}\right) \\
- Z_{\epsilon\gamma\kappa\lambda} & = & \label{e2.30}%
P_{\epsilon\gamma\kappa\lambda} +
2\,A_{\lambda\kappa}A_{\epsilon\gamma}
+
\left(D_{\epsilon\lambda} + A_{\epsilon\lambda}\right)
\left(D_{\gamma\kappa} + A_{\gamma\kappa}\right) -
\left(D_{\gamma\lambda} + A_{\gamma\lambda}\right)
\left(D_{\epsilon\kappa} + A_{\epsilon\kappa}\right)
\end{eqnarray}
}
where
$P_{\epsilon\gamma\kappa\lambda} =
g_{\lambda\sigma}P_{\epsilon\gamma\kappa}^{\hspace*{0.5cm}\sigma} =
- h_{\lambda\sigma}P_{\epsilon\gamma\kappa}^{\hspace*{0.5cm}\sigma}$.
\end{Prop}
\begin{proof}
Relation (\ref{e2.29}) is equal to the relation (3.39) in \cite[p.56]{vlad}, in which a transformation (\ref{e2.21}) of the left-hand side was carried out. From (\ref{e2.29}) we obtain (\ref{e2.30}) by lowering of $\lambda$ by means of $g_{\lambda\sigma}$. The relation $P_{\epsilon\gamma\kappa\lambda} =
- h_{\lambda\sigma}P_{\epsilon\gamma\kappa}^{\hspace*{0.5cm}\sigma}$ is a consequence of (\ref{e2.2}) and (\ref{e2.3}).
\end{proof}

\section{Stationary and static space-times} \label{sec3}%
Now we apply the projection formalism of Section \ref{sec2} to a stationary space-time. Let ${\xi}^{\mu}$ be the timelike Killing field of such a space-time and ${\tau}^{\mu}$ be the timelike unit vector field which is proportional to ${\xi}^{\mu}$, i.e. we have
${\xi}^{\mu} = \phi\,{\tau}^{\mu}$.

\subsection{Proof of Lemma \ref{lemma1.4}}
From ${\xi}^{\mu} = \phi\,{\tau}^{\mu}$ we obtain
\begin{eqnarray} \label{e3.2}%
{\xi}_{\mu ; \nu} & = & ({\partial}_{\nu}\phi){\tau}_{\mu} + \phi{\tau}_{\mu ; \nu}\,.
\end{eqnarray}
and
${\xi}_{[\lambda}{\xi}_{\mu ; \nu]}\;=\;\phi{\tau}_{[\lambda}{\tau}_{\mu}({\partial}_{\nu]}\phi) + {\phi}^2{\tau}_{[\lambda}{\tau}_{\mu ; \nu]}\;=\;{\phi}^2{\tau}_{[\lambda}{\tau}_{\mu ; \nu]}$,
since ${\tau}_{[\lambda}{\tau}_{\mu]} = 0$. Consequently, ${\xi}_{[\lambda}{\xi}_{\mu ; \nu]} = 0$ and ${\tau}_{[\lambda}{\tau}_{\mu ; \nu]} = 0$ are equivalent. $\Box$

\subsection{Proof of Proposition \ref{prop1.2}}
First we show
\begin{Lem}
If ${\tau}^{\mu}$ is proportional to a Killing field ${\xi}^{\mu}$, i.e. ${\xi}^{\mu} = \phi\,{\tau}^{\mu}$ is fulfilled, then
$D_{\mu\nu} = 0$.
\end{Lem}
\begin{proof}
The Killing equation ${\xi}_{(\mu ; \nu)} = 0$ and (\ref{e3.2}) lead to
$0 = {\tau}_{(\mu}{\partial}_{\nu)}\phi + \phi{\tau}_{(\mu ; \nu)}$.
But then (\ref{e2.14}) and (\ref{e2.3}) yield
$D_{\alpha\beta} =
- {\tau}_{(\mu ; \nu)}h_{\alpha}^{\mu}h_{\beta}^{\nu} = {\phi}^{-1}{\tau}_{(\mu}{\partial}_{\nu)}\phi\,h_{\alpha}^{\mu}h_{\beta}^{\nu} = 0$.
\end{proof}

Now Proposition \ref{prop1.2} can be proved in the following way. Because of $D_{\mu\nu} = 0$ we can transform (\ref{e2.30}) into
\begin{eqnarray} \label{e3.5}%
- Z_{\epsilon\gamma\kappa\lambda} & = &
P_{\epsilon\gamma\kappa\lambda} +
2\,A_{\lambda\kappa}A_{\epsilon\gamma} +
A_{\epsilon\lambda}A_{\gamma\kappa} -
A_{\gamma\lambda}A_{\epsilon\kappa}\,.
\end{eqnarray}
If we substitute $A_{\alpha\beta}$ by means of (\ref{e2.16}) in (\ref{e3.5}) and use the notation ${\tau}_{[\mu ; \nu]} = \frac{1}{2}({\tau}_{\mu , \nu} - {\tau}_{\nu , \mu})$, then we obtain (\ref{e1.5}).
We determined the long formula (\ref{e1.5}) by means of the \textsf{Mathematica} package \textsf{Ricci} \cite{ricci3}. The \textsf{Mathematica} notebook of this calculation can be downloaded from \cite{fie21}. $\Box$

\subsection{A remark about synchronized coordinate systems}
\begin{Def}
Let $(M,g)$ be a space-time. Local coordinates $x^0, x^1, x^2, x^3$ are called a \textit{synchronized coordinate system} of $(M,g)$ if the coordinates of $g$ with respect to the $x^{\mu}$ fulfill
\begin{eqnarray} \label{e3.6}%
g_{00}\;=\;1 & \;\;,\;\; & g_{0a}\;=\;0\;\;\;\;(a = 1, 2, 3)\,.
\end{eqnarray}
\end{Def}
Obviously, the basis vector $\tau := {\partial}_0$ of a synchronized coordinate system is a (local) timelike unit vector field.
\begin{Prop} Let $\tau = {\partial}_0$ be the timelike unit vector field of a synchronized coordinate system of a space-time $(M,g)$.
\begin{enumerate}
\item{Then the curl of $\tau$ vanishes, i.e.
${\tau}_{\mu ; \nu} - {\tau}_{\mu ; \nu} =
{\partial}_{\nu}{\tau}_{\mu} - {\partial}_{\mu}{\tau}_{\nu} = 0$.
}
\item{The tensors $F_{\alpha}$ and $A_{\alpha\beta}$ formed from $\tau$ vanish and {\rm (\ref{e2.30})} reduces to
\begin{eqnarray*}
- Z_{\epsilon\gamma\kappa\lambda} & = &
P_{\epsilon\gamma\kappa\lambda} +
D_{\epsilon\lambda}D_{\gamma\kappa} -
D_{\gamma\lambda}D_{\epsilon\kappa}\,.
\end{eqnarray*}
}
\end{enumerate}
\end{Prop}
\begin{proof}
The field $\tau = {\partial}_0$ has the coordinates ${\tau}^0 = 1$, ${\tau}^a = 0$ and ${\tau}_0 = 1$, ${\tau}_a = 0$ with respect to the synchronized coordinate system. (Consider (\ref{e3.6}).) This yields ${\partial}_{\mu}{\tau}_{\nu} = 0$ and the coordinate independent condition
${\tau}_{[\mu ; \nu]} = {\tau}_{[\mu , \nu]} = 0$.
\end{proof}

Around every point $p\in M$ of a space-time $(M,g)$ we can find an infinite set of synchronized coordinate systems. For the $\tau$ of such a synchronized coordinate system we have $F_{\alpha} = 0$, $A_{\alpha\beta} = 0$. However, the set $\mathfrak{T}$ of timelike unit vector fields $\tau$ belonging to a synchronized coordinate system is a proper subset of the set of all timelike unit vector fields of $(M,g)$. For a timelike unit vector field $\tau\not\in\mathfrak{T}$ we have to expect $F_{\alpha} \not= 0$ and/or $A_{\alpha\beta} \not= 0$.


\section{Symmetry classes of tensors} \label{sec4}%
We denote by $\bbK[\calS_r]$ the group ring of the symmetric group $\calS_r$.
\begin{Def}
If $T\in\calT_r V$ and $a = \sum_{p\in\calS_r} a(p)\,p\in\bbK[\calS_r]$, then we denote by $aT$ the $r$-times covariant tensor
\begin{eqnarray} \label{e4.1}%
\hspace*{-0.7cm}(aT)(v_1,\ldots,v_r)\;:=\;\sum_{p\in\calS_r}\,a(p)T(v_{p(1)},\ldots,v_{p(r)})
& \;\;,\;\; &
(aT)_{i_1\ldots i_r}\;=\;\sum_{p\in\calS_r}\,a(p)T_{i_{p(1)}\ldots i_{p(r)}}\,.
\end{eqnarray}
Because of {\rm (\ref{e4.1})}, the group ring elements $a\in\bbK[\calS_r]$ are called \textit{symmetry operators} for the tensors $T\in\calT_r V$.
Further we denote by $\ast:\,\bbK[\calS_r]\rightarrow\bbK[\calS_r]$ the operator
\begin{eqnarray}
\ast:\,a\;=\;\sum_{p\in\calS_r}a(p)p & \mapsto &
a^{\ast}\;:=\;\sum_{p\in\calS_r}a(p)p^{-1}\,.
\end{eqnarray}
\end{Def}
\begin{Def}
Let $\frakr\subseteq\bbK[\calS_r]$ be a right ideal of $\bbK[\calS_r]$. Then the tensor set
\begin{eqnarray}
\calT_{\frakr} & := & \{ aT\;|\;a\in\frakr\;,\;T\in\calT_r V\}
\end{eqnarray}
is called the \textit{symmetry class} of $r$-times covariant tensors defined by $\frakr$. $\calT_{\frakr}$ is called irreducible iff $\frakr$ is minimal.
\end{Def}
\begin{Prop}\footnote{See \cite[pp.127]{boerner}, \cite[Lemma III.2.2]{fie16}, \cite{fie18}.}
Let $e\in\bbK[\calS_r]$ be a generating idempotent of a right ideal $\frakr\subseteq\bbK[\calS_r]$. Then a tensor $T\in\calT_r V$ lies in the symmetry class $\calT_{\frakr}$ of $\frakr$ iff $eT = T$.
\end{Prop}
Important special symmetry operators are Young symmetrizers, which are defined by means of Young tableaux.
A {\itshape Young tableau} $t$ of $r\in\bbN$ is an arrangement of $r$ boxes such that
\begin{enumerate}
\item{the numbers ${\lambda}_i$ of boxes in the rows $i = 1 , \ldots , l$ form a decreasing sequence
${\lambda}_1 \ge {\lambda}_2 \ge \ldots \ge {\lambda}_l > 0$ with
${\lambda}_1 + \ldots + {\lambda}_l = r$,}
\item{the boxes are filled with the numbers $1, 2, \ldots , r$ in any order.}
\end{enumerate}
For instance, the following graphic shows a Young tableau of $r = 15$.
{\scriptsize
\[\left.
\begin{array}{cc|c|c|c|c|c|c}
\cline{3-7}
{\lambda}_1 = 5 & \;\;\; & 11 & 2 & 5 & 4 & 12 & \\
\cline{3-7}
{\lambda}_2 = 4 & \;\;\; & 9 & 6 & 13 & 3 \\
\cline{3-6}
{\lambda}_3 = 4 & \;\;\; & 8 & 15 & 1 & 7 \\
\cline{3-6}
{\lambda}_4 = 2 & \;\;\; & 10 & 14 \\
\cline{3-4}
\end{array}\right\}\;=\;t\,.
\]
}
Obviously, the unfilled arrangement of boxes, the {\itshape Young frame}, is characterized by a partition
$\lambda = ({\lambda}_1 , \ldots , {\lambda}_l) \vdash r$ of $r$.

If a Young tableau $t$ of a partition $\lambda \vdash r$ is given, then the {\itshape Young symmetrizer} $y_t$ of $t$ is defined by\footnote{We use the convention $(p \circ q) (i) := p(q(i))$ for the product of two permutations $p, q$.}
\begin{eqnarray}
y_t & := & \sum_{p \in {\calH}_t} \sum_{q \in {\calV}_t} \mathrm{sign}(q)\, p \circ q
\end{eqnarray}
where ${\calH}_t$, ${\calV}_t$ are the groups of the {\itshape horizontal} or
{\itshape vertical permutations} of $t$ which only permute numbers within rows or columns of $t$, respectively. The Young symmetrizers of $\bbK [{\calS}_r]$ are essentially idempotent and define decompositions
\begin{eqnarray}
\bbK [{\calS}_r] \;=\;
\bigoplus_{\lambda \vdash r} \bigoplus_{t \in {\calS\calT}_{\lambda}}
\bbK [{\calS}_r]\cdot y_t
& \;\;,\;\; &
\bbK [{\calS}_r] \;=\;
\bigoplus_{\lambda \vdash r} \bigoplus_{t \in {\calS\calT}_{\lambda}}
y_t \cdot \bbK [{\calS}_r] \label{e4.5}
\end{eqnarray}
of $\bbK [{\calS}_r]$ into minimal left or right ideals. In (\ref{e4.5}), the symbol ${\calS\calT}_{\lambda}$ denotes the set of all standard tableaux of the partition $\lambda$. \textit{Standard tableaux} are Young tableaux in which the entries of every row and every column form an increasing number sequence.\footnote{About Young symmetrizers and
Young tableaux see for instance the references given in \cite[Footnote 6]{fie04a}.}


The inner sums of (\ref{e4.5}) are minimal two-sided ideals
\begin{eqnarray}
{\mathfrak a}_{\lambda} & := &
\bigoplus_{t \in {\calS\calT}_{\lambda}}
\bbK [{\calS}_r]\cdot y_t
\;=\;
\bigoplus_{t \in {\calS\calT}_{\lambda}}
y_t \cdot \bbK [{\calS}_r]
\end{eqnarray}
of $\bbK [{\calS}_r]$.
The set of all Young symmetrizers $y_t$ which lie in ${\mathfrak a}_{\lambda}$ is equal to the set of all $y_t$ whose tableau $t$ has the frame $\lambda\vdash r$.
Furthermore two minimal left ideals $\frakl_1, \frakl_2\subseteq\bbK[\calS_r]$ or two minimal right ideals $\frakr_1, \frakr_2\subseteq\bbK[\calS_r]$ are \textit{equivalent} iff they lie in the same ideal ${\mathfrak a}_{\lambda}$.
Now we say that a symmetry class $\calT_{\frakr}$ \textit{belongs to} $\lambda\vdash r$ iff
$\frakr\subseteq{\mathfrak a}_{\lambda}$.

S.A. Fulling, R.C. King, B.G.Wybourne and C.J. Cummins showed in \cite{full4}
\begin{Thm} {\rm\bf (Fulling, King, Wybourne, Cummins)}\footnote{Theorem \ref{thm4.3} is a special case of a more general theorem about symmetry classes of $R_{i j k l\,;\,(s_1 \ldots s_u)}$ in \cite{full4}.} \label{thm4.3} \\
Let $y_t$ be the Young symmetrizer of the standard tableau {\rm (\ref{e4.8})}.
Then a tensor $T\in\calT_4 V$ lies in
$\calA(V)$ if and only if
$\textstyle{\frac{1}{12}}\,y_t^{\ast} T\;=\;T$.
\end{Thm}
The group ring element $\frac{1}{12}\,y_t^{\ast}$ is a primitive idempotent. Since $y_t\in {\mathfrak a}_{(2^2)}$, the minimal ideals $\frakl := \bbK[\calS_4]\cdot y_t$, $\frakr = \frakl^{\ast} = y_t^{\ast}\cdot\bbK[\calS_4]$ satisfy $\frakl , \frakr\subset\mathfrak{a}_{(2^2)}$, i.e. the symmetry class $\calA(V)$ of algebraic curvature tensors belongs to the partition $\lambda = (2^2)\vdash 4$.

The following proposition guarantees that we can use Littlewood-Richardson products to determine information about symmetry classes which contain product tensors such as ${\tau}_{\lambda}{\tau}_{[\mu ; \nu]}$ and $F_{\kappa}{\tau}_{\lambda}{\tau}_{[\mu ; \nu]}$.
\begin{Prop}\hspace{- 1mm}\footnote{See B. Fiedler \cite[Sec.III.3.2]{fie16} and B. Fiedler \cite{fie17}.} \label{prop2.8}%
Let ${\frakr}_i \subseteq {\bbK} [{\calS}_{r_i}]$ $(i = 1,\ldots , m)$ be right ideals
and
$T^{(i)} \in {\calT}_{{\frakr}_i^{\ast}} \subseteq {\calT}_{r_i} V$ be $r_i$-times 
covariant tensors from the symmetry classes characterized by the
${\frakr}_i$.
Consider the product
\begin{eqnarray}
T & := & T^{(1)} \otimes\ldots\otimes T^{(m)} \;\in\; {\calT}_r V
\;\;\;,\;\;\;
r := r_1 + \ldots + r_m \,. \label{equ2.17}%
\end{eqnarray}
For every $i$ we define an embedding
\begin{eqnarray}
{\iota}_i : {\calS}_{r_i} \rightarrow {\calS}_r
& \;\;,\;\; &
({\iota}_i s)(k) := 
\left\{
\begin{array}{ll}
{\Delta}_i + s(k - {\Delta}_i) & {\rm if}\;\; r_{i-1} < k \le r_i \\
k & {\rm else}
\end{array}
\right.
\end{eqnarray}
where ${\Delta}_i := r_0 + \ldots + r_{i-1}$ and $r_0 := 0$.
Then all product tensors {\rm (\ref{equ2.17})} belong to the symmetry class $\calT_{\frakr}$ of
the right ideal
\begin{eqnarray}
{\frakr} & := &
\left({\bbK} [{\calS}_r]\cdot{\calL}\bigl\{{\tilde{\frakl}}_1 \cdot\ldots\cdot {\tilde{\frakl}}_m
\bigr\}\right)^{\ast} \; = \;
\left({\bbK} [{\calS}_r]\cdot\bigl({\tilde{\frakl}}_1 \otimes\ldots\otimes
{\tilde{\frakl}}_m \bigr)\right)^{\ast} \label{equ2.19}%
\end{eqnarray}
where ${\tilde{\frakl}}_i := {\iota}_i ({\frakl}_i)$ are the embeddings of
the left ideals ${\frakl}_i = {\frakr}_i^{\ast}$ into ${\bbK} [{\calS}_r]$ induced by the ${\iota}_i$. If
$\dim V \ge r$, then the right ideal $\frakr$ does not contain a proper right subideal $\tilde{\frakr}\subset\frakr$ such that all tensors {\rm (\ref{equ2.17})} lie also in the symmetry class $\calT_{\tilde{\frakr}}$ of $\tilde{\frakr}$.
\end{Prop}
Let $\regrep_G : G \rightarrow GL({\bbK}[G])$ denote the {\it regular
representation} of a finite group $G$ defined by
$\regrep_g (f) := g \cdot f$, $g \in G$, $f \in {\bbK}[G]$. If we use the 
left ideals ${\frakl}_i = {\frakr}_i^{\ast}$, ${\frakl} = {\frakr}^{\ast}$
to define subrepresentations
${\alpha}_i := \regrep_{{\calS}_{r_i}} |_{{\frakl}_i}$,
$\beta := \regrep_{{\calS}_r} |_{\frakl}$,
then the representation $\beta$ is equivalent to the {\itshape Littlewood-Richardson 
product} of the ${\alpha}_i$ (see B. Fiedler \cite[Sec.III.3.2]{fie16}):
\begin{eqnarray}
\beta & \sim & {\alpha}_1 {\alpha}_2\ldots {\alpha}_m\;:=\;
{\alpha}_1 \,\#\ldots\#\, {\alpha}_m \uparrow {\calS}_r \label{equ2.20}
\end{eqnarray}
('$\#$' denotes the outer tensor product of the above representations.) This result corresponds to statements of S.A. Fulling et al. \cite{full4}.
Relation (\ref{equ2.20}) allows us to determine information about the structure of the right ideal (\ref{equ2.19}) by means of the {\itshape Littlewood-Richardson rule}\footnote{About the Littlewood-Richardson rule see for instance the references in \cite[Footnote 13]{fie04a}.}.

\section{On symmetry classes for ${\tau}_{\lambda}{\tau}_{[\mu ; \nu]}$ and
$F_{\kappa}{\tau}_{\lambda}{\tau}_{[\mu ; \nu]}$}
\subsection{Proof of Theorem \ref{thm1.6}}
Because of Proposition \ref{prop2.8} we can determine information about a symmetry class containing ${\tau}_{\lambda}{\tau}_{[\mu ; \nu]}$ from the corresponding Littlewood-Richardson product. Since the symmetries of ${\tau}_{\lambda}$ and ${\tau}_{[\mu ; \nu]}$ are generated by the Young symmetrizers of the tableaux
{\scriptsize
$\begin{array}{|c|}
\hline
 1 \\
\hline
\end{array}
\;,\;
\begin{array}{|c|}
\hline
 1 \\
\hline
 2 \\
\hline
\end{array}$
}
we have to calculate the Littlewood-Richardson product $[1^2][1]$. The use of the Littlewood-Richardson rule yields the graphic
{\scriptsize
\begin{eqnarray*}
\begin{array}{|c|}
\hline
  \\
\hline
 \\
\hline
\end{array}
\;\;\times\;
\begin{array}{|c|}
\hline
  \\
\hline
\end{array}
& \;\;\;\sim\;\;\; &
\begin{array}{|c|c|c}
\cline{1-2}
  &  & \\
\cline{1-2}
 \\
\cline{1-1}
\end{array}
+\;\;\;
\begin{array}{|c|}
\hline
 \\
\hline
 \\
\hline
 \\
\hline
\end{array}
\end{eqnarray*}
}
which can be translated into (\ref{e1.6}). The relation (\ref{e1.6}) means that the right ideal $\frakr$ of the symmetry class of ${\tau}_{\lambda}{\tau}_{[\mu ; \nu]}$ possesses a decomposition $\frakr = \frakr_1\oplus\frakr_2$ into 2 minimal right ideals which belong to the partitions $(2\,1), (1^3)\vdash 3$.

Now we see that information about the right ideal $\frakr$ of a symmetry class for $F_{\kappa}{\tau}_{\lambda}{\tau}_{[\mu ; \nu]}$ can be gained from the Littewood-Richardson products $[2\,1][1]$ and $[1^3][1]$. For these products the Littlewood-Richardson rule yields (\ref{e1.14}) since it leads to the graphics
{\scriptsize
\begin{eqnarray*}
\begin{array}{|c|c|c}
\cline{1-2}
 & & \\
\cline{1-2}
 \\
\cline{1-1}
\end{array}
\hspace{-0.2cm}\times\;
\begin{array}{|c|}
\hline
 \\
\hline
\end{array}
\;\;\;\sim\;\;\;
\begin{array}{|c|c|c|c}
\cline{1-3}
 & & & \\
\cline{1-3}
 \\
\cline{1-1}
\end{array}
\hspace{-0.1cm}+\;\;
\begin{array}{|c|c|}
\hline
 & \\
\hline
 & \\
\hline
\end{array}
\;\;\;+\;\;
\begin{array}{|c|c|c}
\cline{1-2}
 & & \\
\cline{1-2}
 \\
\cline{1-1}
 \\
\cline{1-1}
\end{array}
& \hspace*{0.5cm},\hspace{1cm} &
\begin{array}{|c|}
\hline
 \\
\hline
 \\
\hline
 \\
\hline
\end{array}
\hspace{0.2cm}\times\;
\begin{array}{|c|}
\hline
 \\
\hline
\end{array}
\;\;\;\sim\;\;\;
\begin{array}{|c|c|c}
\cline{1-2}
 & & \\
\cline{1-2}
 \\
\cline{1-1}
 \\
\cline{1-1}
\end{array}
+\;\;\;\;
\begin{array}{|c|}
\hline
 \\
\hline
 \\
\hline
 \\
\hline
 \\
\hline
\end{array}\;.
\end{eqnarray*}
}

Theorem \ref{thm4.3} tells us that the symmetry class $\calA(V)$ is defined by the right ideal $\frakr = y_t^{\ast}\cdot\bbK[\calS_4]$ generated by the Young symmetrizer of the tableau (\ref{e4.8}). This right ideal $\frakr$ is minimal and belongs to $(2^2)\vdash 4$ because $y_t\in\mathfrak{a}_{(2^2)}$. But only $[2\,1][1]$ possesses a part $[2^2]$ which belongs to a minimal right ideal of $(2^2)$. This proves the last assertion of Theorem \ref{thm1.6}. $\Box$

\subsection{Proof of Theorem \ref{thm1.8}}
One and only one symmetry class $\calT_{\tilde{\frakr}}$ of $\calT_3 V$ belongs to the partition $(1^3)\vdash 3$. It is defined by the right ideal
$\tilde{\frakr} := y_{\tilde{t}}\cdot\bbK[\calS_3]$ where $y_{\tilde{t}}$ is the Young symmetrizer of the Young tableau
{\scriptsize
\begin{eqnarray*}
\tilde{t} & = &
\begin{array}{|c|}
\hline
1 \\
\hline
2 \\
\hline
3 \\
\hline
\end{array}\,.
\end{eqnarray*}
}
$\frac{1}{6}y_{\tilde{t}}$ yields the alternation of a tensor from $\calT_3 V$, i.e. it holds
$\frac{1}{6}(y_{\tilde{t}}T)_{\lambda\mu\nu} = {\tau}_{[\lambda}{\tau}_{\mu ; \nu]}$ if we abbreviate $T_{\lambda\mu\nu} := {\tau}_{\lambda}{\tau}_{[\mu ; \nu]}$. Thus, ${\tau}_{[\lambda}{\tau}_{\mu ; \nu]} = 0$ means that the part of ${\tau}_{\lambda}{\tau}_{[\mu ; \nu]}$ which lies in $\calT_{\tilde{\frakr}}$ vanishes. $\Box$

\subsection{Proof of Theorem \ref{thm1.9}}
Because of (\ref{e1.14}) the Littlewood-Richardson product $[1^3][1]$ has no part that belongs to the partition $(2^2)\vdash 4$. This leads us to the assumtion that the terms ${\tau}_{[\lambda}{\tau}_{\mu ; \nu]}$ will possibly fall out of (\ref{e1.5}) if we carry out the substitution (\ref{e1.17}) in (\ref{e1.5}). We verified this by a \textsf{Mathematica} calculation using the tensor package \textsf{Ricci} \cite{ricci3}. In this calculation we showed that the lines 3, 4, 5 of (\ref{e1.5}) vanish if we replace all expressions ${\tau}_{\lambda}{\tau}_{[\mu ; \nu]}$ by the coordinates of an arbitrary alternating tensor $a_{\lambda\mu\nu}$ of order 3. A record of this calculation can be found in \cite[curvterms.nb]{fie21}. $\Box$

\subsection{Proof of Theorem \ref{thm1.10}}
The proof uses results from \cite{fie03b,fie04a}. In \cite{fie03b} we showed by means of \textit{discrete Fourier transforms} for symmetric groups that every minimal right ideal $\frakr\subset\mathfrak{a}_{(2\,1)}\subset\bbK[\calS_3]$ is generated by exactly one element of the following set of (primitive) idempotents\footnote{In \cite{fie03b,fie04a} the idempotents ${\zeta}_{\nu}$ were denoted by ${\xi}_{\nu}$.}
\begin{eqnarray*}
{\zeta}_{\nu} & := & \frac{1}{3}\left\{[1,2,3] + \nu [1,3,2] + (1 - \nu)[2,1,3] \right.\\
 & & \left.- \nu [2,3,1] + (-1 + \nu)[3,1,2] - [3,2,1]\right\}\;\;,\;\;\nu\in\bbK \\
\eta & := & \frac{1}{3}\left\{[1,2,3] - [2,1,3] - [2,3,1] + [3,2,1]\right\}\,.
\end{eqnarray*}
On the other hand
${\theta}_{\lambda\mu\nu} = {\tau}_{\lambda}{\tau}_{[\mu ; \nu]} -
{\tau}_{[\lambda}{\tau}_{\mu ; \nu]}$ is generated by the idempotent symmetry operator
$\rho := \frac{1}{2}\left\{[1,2,3] - [1,3,2]\right\} -
\frac{1}{6}\sum_{p\in\calS_3}\,\mathrm{sign}(p)\,p$
from ${\tau}_{\lambda}{\tau}_{\mu ; \nu}$.
We showed by means of the \textsf{Mathematica} package \textsf{PERMS} \cite{fie10} in the notebook \cite[curvterms.nb]{fie21} that
$\eta\cdot\rho \not= \rho$, $\rho\cdot\eta \not= \eta$ and
\begin{eqnarray*}
{\zeta}_{\nu}\cdot\rho = \rho\;\;\mathrm{and}\;\;\rho\cdot{\zeta}_{\nu} = {\zeta}_{\nu} & \Leftrightarrow & \nu = -1\,.
\end{eqnarray*}
Consequently, the tensor ${\theta}_{\lambda\mu\nu}$ belongs to the symmetry class $\calT_{\frakr}$ which is defined by the right ideal
$\frakr := {\zeta}_{-1}\cdot\bbK[\calS_3]$ with generating idempotent ${\zeta}_{-1}$.

This right ideal $\frakr$ is different from the right ideal $\frakr_0$ in Theorem \ref{thm1.3} since a result of \cite{fie04a} says that $\frakr_0 = {\zeta}_{1/2}\cdot\bbK[\calS_3]$.
In the notebook \cite[part16a.nb]{fie21} for the paper \cite{fie04a} we proved that all tensors of the symmetry class of ${\zeta}_{-1}$ satisfy the identities (\ref{e1.19}). Using (\ref{e1.19}) we transformed (\ref{e1.18}) into (\ref{e1.20}) by means of \textsf{Ricci} \cite{ricci3} in \cite[curvterms.nb]{fie21}. Finally, we see by a comparison with formulas in \cite[Sec.4.3]{fie04a} that the last line of (\ref{e1.20}) is equal to $\frac{1}{2} y_t^{\ast}(\theta\otimes F)_{\kappa\lambda\mu\nu}$, where $t$ is the tableau (\ref{e4.8}). $\Box$\\*[0.1cm]
\noindent {\bf Remark :} In \cite{fie04a} we showed that the symmetry operator ${\zeta}_{-1}$ produces tensors $U_{\lambda\mu\nu}$ of a $(2\,1)$-symmetry class which admits the \textit{index commutation symmetry} $U_{\lambda\nu\mu} = - U_{\lambda\mu\nu}$. Obviously the above tensor ${\theta}_{\lambda\mu\nu}$ possesses this index commutation symmetry. A further result of \cite{fie04a} says that then the coordinates of algebraic curvature tensors
$y_t^{\ast}(w\otimes U)$ ($t$ given by (\ref{e4.8})) can be reduced to a sum with a minimal length of 4 summands. This effect can be observed in the above reduction of the number of terms $F_{\kappa}{\theta}_{\lambda\mu\nu}$ in formula (\ref{e1.18}), too.\\*[0.3cm]
\noindent {\bf Acknowledgements.} I would like to thank Prof. P. B. Gilkey for important and helpful discussions and for valuable suggestions for future investigations.

\section*{References}

\end{document}